%% file: Preprint_4_QuasiMonoticPresFin.tex
\documentclass[12pt]{article}

\usepackage{erics_preprints}

\input xy
\xyoption{all}
\xyoption{poly}
\xyoption{matrix}

\begin{document}

\title{An Alternative Presentation of the Symmetric-Simplicial Category}
\author{Eric Ram\'on Antokoletz}
\date{}
\maketitle

\begin{abstract}

The category $\fin$ of {\it symmetric-simplicial operators} is obtained
by enlarging the category $\ord$ of monotonic functions between the sets
$\{0,1,\dots n\}$ to include all functions between the same sets. Marco
Grandis \cite{Grandis} has given a presentation of $\fin$ using the standard
generators $d_i$ and $s_i$ of $\ord$ as well
as the adjacent transpositions $t_i$ which generate the permutations in
$\fin$. The purpose of this note is to establish an alternative
presentation of $\fin$ in which the codegeneracies $s_i$ are replaced by
{\it quasi-codegeneracies} $u_i$. 
We
also prove a unique factorization theorem for products of $d_i$ and
$u_j$ analogous to the standard unique factorizations in $\ord$.  This presentation has been used by the
author to construct {\it symmetric
hypercrossed complexes} (to be published elsewhere) which are algebraic models for homotopy types of
spaces based on the hypercrossed complexes of \cite{CarrascoCegarra}. 

\end{abstract}

\tableofcontents
		
\setcounter{secnumdepth}{2}
\setcounter{section}{0}
	
	\section{Introduction}						\input{IntroFinAltPres2}

	\section{Grandis's Presentation \\ of the Symmetric-Simplicial Category}	\input{SymmMooreComp}

	\section{Proof of the Alternative Presentation \\ of the Symmetric-Simplicial Category}\input{SMCProof}
	\section{The Algebra of the Symmetric-Simplicial Category}	\input{SymAlgMultIndex}

\bibliographystyle{amsalpha}
\bibliography{Preprint_4_QuasiMonoticPresFin}

\end{document}

%% file: IntroFinAltPres2.tex


	\label{IntroFinAltPres}

In order to motivate the subject of this note, we bring together two
distinct lines of historical development. First, recall that $\ord$ is
the category whose objects are the standard finite ordered sets
	$$[n] := \{0,1,\ldots, n\} \mfor n \ge 0$$
and whose morphisms are all monotonic functions
	$$f: [n] \lra [m]$$
	$$i < j \implies f(i) \le f(j).$$
Simplicial sets, which are by definition contravariant functors from
$\ord$ to the category $\set$ of all sets and mappings, are used in
Homotopy Theory and related fields as combinatorial models for
topological spaces, among other things. For this reason, $\ord$ is
often referred to as {\it the simplicial category} and denoted $\Delta$,
and the category of simplicial sets is then denoted
$\set_\Delta:=\fun(\Delta^{op},\set)$.

In Pursuing Stacks (\cite{Groth}), Alexander Grothendieck proposed replacing $\ord$ in the
definition of simplicial set with an arbitrary small category $\Gamma$
and looking for Quillen model category structures on the category
$\set_\Gamma$ of $\Gamma$-sets (defined as $\set_\Gamma:=\gammaset$) in
order to investigate the possibilities for doing homotopy theory there.
He laid special emphasis on certain geometrically motivated examples of
$\Gamma$, including the category which is the subject of this paper,
namely the category denoted $\fin$ (denoted by him $\wt\Delta$) whose
objects are the same as those of $\ord$ but whose morphisms consist of
{\it all} functions $f: [n] \lra [m]$ for each $m,n\ge 0$. 

A short while later, W.G.Dwyer, Michael Hopkins and Daniel Kan proved a result 
showing that for a certain class of categories $\Gamma$, including
$\Gamma = \fin$, one may define a model structure on $\set_\Gamma$ such
that the resulting homotopy theory is equivalent to the usual one on
$\set_\Delta$ (see \cite{DwyerHopkinsKan}). A later observation of F.William Lawvere in
\cite{Lawvere} also suggested studying $\set_\Gamma$ for $\Gamma =
\fin$, inspiring Marco Grandis to take up the subject (\cite{Grandis2},
\cite{Grandis}, \cite{Grandis3}, \cite{Grandis1}, \cite{GrandisNew}).
Note that $\fin$ contains the group $\sym[n]$ of all permutations of the set $[n]$
for each $n\ge0$, as well as the category $\ord$. For this reason,
$\fin$ is referred to by Grandis as the {\it symmetric-simplicial
category}, and we shall do so as well.  

We turn briefly to the other line of historical development relevant for us 
here. Motivated by the fact that the category $\sgrp := \fun(\Delta^{op},\grp)$
of simplicial groups possesses a homotopy theory equivalent to that of
pointed connected spaces (they play the role of loop spaces, see
\cite{May} or \cite{GoerssJardine}), P. Carrasco and A. M. Cegarra
discovered a nonabelian Dold-Kan theorem for simplicial groups and used
it to describe homological-algebraic models for classical homotopy types
(\cite{CarrascoCegarra}), which they dubbed {\it hypercrossed
complexes}. 

Since the author has shown (to appear elsewhere) that
symmetric-simplicial groups also have a homotopy theory equivalent to
that of pointed connected spaces, it is of interest to ask what sort of
homological-algebraic objects can arise from nonabelian Dold-Kan
decompositions (in the sense of \cite{CarrascoCegarra}) of
symmetric-simplicial groups. The author has shown (\cite{preprint2})
that, in addition to the decompositions obtained via a direct
application of \cite{CarrascoCegarra} to the underlying simplicial group
of a symmetric-simplicial group, there also exist new Dold-Kan
decompositions which can be obtained by making judicious use of the
algebra of the category $\fin$. These decompositions give rise in turn
to new homological-algebraic models for homotopy types, which we call
{\it symmetric hypercrossed complexes}, that are simpler than the
original hypercrossed complexes in the sense that a great deal of the
algebraic data constituting them vanishes (to appear elsewhere).

The new Dold-Kan decompositions are obtained using an
alternative presentation of the category $\fin$, whose verification is
the main purpose of the present note. Grandis gave a presentation in
\cite{Grandis} (reviewed in section \ref{SymmMooreComp} below) of $\fin$
that uses the standard presentation of $\ord$ as well as the Moore
presentations of the symmetric groups $\sym[n]$ via adjacent
transpositions. In the alternative presentation of $\fin$, the monotonic
elementary codegeneracies $s_i\in\ord$ in Grandis's presentation are
replaced by certain nonmonotonic surjections $u_i\in\fin$ which we call
the {\it elementary quasi-codegeneracies} (see Definition
\ref{quasicodegens}). 

In section \ref{SMCProof} we prove this alternative presentation by
relating it directly to Grandis's presentation. In section
\ref{SymAlgMultIndex}, we also show that the morphisms of the
subcategory of $\fin$ generated by the $d_i$ and $u_i$ are characterized
by the following two conditions.
	\begin{itemize}
		\item{They take $0$ to $0$.}
		\item{They are strictly monotonic outside of the preimage of $0$.}
	\end{itemize}
We call such morphisms {\it
quasi-monotonic} and denote the subcategory of $\fin$ consisting of
quasi-monotonic functions by ${\bf qOrd}$. Finally we show that $\bf
qOrd$ admits unique factorizations analogous to those of $\ord$. These
results are relied upon in \cite{preprint2} to derive the alternative
Dold-Kan decompositions for symmetric-simplicial groups mentioned above.

%% file: SymmMooreComp.tex

\label{SymmMooreComp}

We begin by recalling the well-known presentation of category $\ord$ via
generators and relations (see \cite{May}, \cite{Lamotke},
\cite{MacLane}, \cite{GoerssJardine} et. al.). The generators are given
in the following definition.

\defn{ The {\it elementary coface} maps are defined by
	$$d_i = d_i^{(n)}: [n-1] \lra [n] ~ \mfor n \ge 1 \mand 0 \le i \le n$$
	$$k \mapsto \left\{\begin{array} {cl}
				k & \mfor k \le i \\
				k+1 & \mfor k > i
			\end{array}\right.$$
	and the {\it elementary codegeneracy} maps are defined as follows.
	$$s_i = s_i^{(n)}: [n+1] \lra [n] ~ \mfor n \ge 0 \mand 0 \le i \le n$$
	$$k \mapsto \left\{\begin{array} {cl}
				k & \mfor k \le i \\
				k-1 & \mfor k > i
			\end{array}\right.$$}

\rmk{ One may put this definition into words by saying that $d_i$ is the
unique monotonic injection $[n-1]\rightarrow[n]$ whose image contains
everything except the element $i$, and $s_i$ is the unique monotonic
surjection $[n+1]\rightarrow[n]$ for which each range element has a
single pre-image except for the element $i$, which has two pre-images.}

These generators satisfy the following {\it cosimplicial identities}. It
is proved in \cite{MacLane} that $\ord$ is isomorphic to the abstract
category obtained by imposing the cosimplicial identities on the free category
having objects $[n]$ and generators $s_i^{(n)}, d_i^{(n)}$. 

	$$\text{{\bf \underline{The Cosimplicial Identities.}}}$$
	$$\begin{array}{lcl}
		{d_id_j = \left\{\begin{array}{cl}
				d_{j+1}d_i & ~~\textrm{if}~~ i \le j \\
				d_{j}d_{i-1} & ~~\textrm{if}~~ i > j 
			\end{array}\right.} &&
		{d_is_j = \left\{\begin{array}{cl}
				s_{j+1}d_i & ~~\textrm{if}~~ i \le j \\
				s_jd_{i+1} & ~~\textrm{if}~~ i > j
			\end{array}\right.} \\[10pt]  \\[-10pt]
		{s_id_j = \left\{\begin{array}{cl}
				d_{j-1}s_i & ~~\textrm{if}~~  i \le j-2 \\
				\text{id} & ~~\textrm{if}~~ i = j-1 \mathrm{~or~} j \\
				d_{j}s_{i-1} & ~~\textrm{if}~~ i \ge j+1 
			\end{array}\right.} &&
		{s_is_j = \left\{\begin{array}{cl}
				s_{j-1}s_i & ~~\textrm{if}~~  i < j \\
				s_{j}s_{i+1} & ~~\textrm{if}~~ i \ge j
			\end{array}\right.}
	\end{array}\\[10pt]$$

\rmk{
	These identities are usually written in a nonredundant form. Here, and
	in all other presentations below, we have included all possible
	situations that arise when interchanging two generators, thus incurring
	a certain amount of redundancy.
}

Traditionally, the action of $\ord$ on simplicial objects is written on
the {\it left}. This necessitates reversing the cosimplicial identities
given above (and then reorganizing indices). The reversed identities are
called the {\it simplicial identities} and are included here for
reference.
	$$\text{{\bf \underline{The Simplicial Identities.}}}$$
	$$\begin{array}{ll}
		{d_id_j = \left\{\begin{array}{cl}
				d_{j-1}d_{i} & ~~\textrm{if}~~ i < j \\
				d_jd_{i+1} & ~~\textrm{if}~~ i \ge j  
			\end{array}\right.} &
		  \\[12pt]
		{s_is_j = \left\{\begin{array}{cl}
				s_{j+1}s_{i} & ~~\textrm{if}~~ i \le j \\
				s_js_{i-1} & ~~\textrm{if}~~  i > j 
			\end{array}\right.} & \\[13pt]
		{d_is_j = \left\{\begin{array}{cl}
				s_{j-1}d_{i} & ~~\textrm{if}~~ i < j \\
				\text{id} & ~~\textrm{if}~~ i = j \mathrm{~or~} j+1 \\
				s_jd_{i-1} & ~~\textrm{if}~~  i \ge j+2
			\end{array}\right.} &
			s_id_j = \left\{\begin{array}{cl}
				d_{j+1}s_i & ~~\textrm{if}~~ i < j \\
				d_js_{i+1} & ~~\textrm{if}~~ i \ge j 
			\end{array}\right. \\
\end{array}$$

\rmk{
	The table above (as well as the others to come) is arranged so that all
	identities in the right column follow from the identities to their left.
	Some redundancies also remain within the left column. 
}

A presentation of $\fin$ via generators and relations has been given by
Marco Grandis in \cite{Grandis}. In addition to the generators $d_i$ and
$s_i$ of $\ord$, his presentation also makes use of the following
generating permutations.

\defn{ The {\it adjacent transpositions} are defined as follows.
	$$t_i = t_i^{(n)}: [n] \lra [n] ~ \mfor n \ge 1 \mand 0 \le i \le n-1$$
	$$k \mapsto \left\{\begin{array} {cl}
				k & \mfor k \neq i,~i+1 \\
				i+1 & \mfor k = i \\
				i & \mfor k = i+1
			\end{array}\right.$$
}

These transpositions satisfy certain relations constituting a well-known
presentation of the symmetric group on $n+1$ elements, ascribed to the
American mathematician E.H.Moore (1862-1932).
	$$\begin{array}{rcl}
			t_i^2 & = &  \text{id} \\
			t_it_j & = & t_jt_i \text{~~if~~} \abs{i-j} \ge 2 \\
			t_it_{i+1}t_i & = & t_{i+1}t_it_{i+1}
	\end{array}$$
In addition to these as well as the simplicial identities, Grandis's presentation also includes
relations allowing one to interchange a transposition with a face or
degeneracy operator.  His relations are given below in contravariant form,
that is, as the relations defining $\fin^{op}$, so that they are
suitable for writing the action on a symmetric-simplicial object on the {\it left}. 
\pagebreak
$$\text{{\bf \underline{The Symmetric-Simplicial Identities (Grandis).}}}$$
	$$\begin{array}{ll}
		{d_id_j = \left\{\begin{array}{cl}
				d_{j-1}d_{i} & ~~\textrm{if}~~ i < j \\
				d_jd_{i+1} & ~~\textrm{if}~~ i \ge j  
			\end{array}\right.} & \\[15pt]
		{d_is_j = \left\{\begin{array}{cl}
				s_{j-1}d_{i} & ~~\textrm{if}~~ i < j \\
				\text{id} & ~~\textrm{if}~~ i = j \mathrm{~or~} j+1 \\
				s_jd_{i-1} & ~~\textrm{if}~~  i \ge j+2
			\end{array}\right.} &
		{s_id_j = \left\{\begin{array}{cl}
				d_{j+1}s_i & ~~\textrm{if}~~ i < j \\
				d_js_{i+1} & ~~\textrm{if}~~ i \ge j 
			\end{array}\right.} \\[20pt]
		{s_is_j = \left\{\begin{array}{cl}
				s_{j+1}s_{i} & ~~\textrm{if}~~ i \le j \\
				s_js_{i-1} & ~~\textrm{if}~~  i > j 
			\end{array}\right.} \\[13pt]
		{t_it_j = \left\{\begin{array}{cl}
			\id & ~~\textrm{if}~~ i = j \\
			t_jt_i & ~~\textrm{if}~~ \abs{i-j} \ge 2 \\
			(t_jt_i)^2 & ~~\textrm{if}~~ \abs{i-j} = 1
			\end{array}\right.} \\[20pt]
		{d_it_j = \left\{\begin{array}{cl}
			t_{j-1}d_{i} & ~~\textrm{if}~~ i < j \\
			d_{i+1} & ~~\textrm{if}~~ i = j \\
			d_{i-1} & ~~\textrm{if}~~ i = j+1 \\
			t_jd_i & ~~\textrm{if}~~ i \ge j+2
			\end{array}\right.}  &
		{t_id_j = \left\{\begin{array}{cl}
			d_jt_i & ~~\textrm{if}~~ i \le j-2 \\
			d_jt_{i+1}t_it_{i+1} & ~~\textrm{if}~~ i = j-1 \\
			d_jt_{i+1} & ~~\textrm{if}~~ i \ge j
			\end{array}\right.} \\[27pt]
		 {t_is_j = \left\{\begin{array}{cl}
			s_jt_i & ~~\textrm{if}~~ i \le j-2 \\
			t_{i+1}s_{j-1}t_i & ~~\textrm{if}~~ i = j-1 \\
			s_j & ~~\textrm{if}~~ i = j \\
			t_{i-1}s_{j+1}t_{i-1} & ~~\textrm{if}~~ i = j+1 \\
			s_jt_{i-1} & ~~\textrm{if}~~ i \ge j+2
			\end{array}\right.} &
		{s_it_j = \left\{\begin{array}{cl}
			t_{j+1}s_{i} & ~~\textrm{if}~~ i < j \\
			t_{j+1}t_js_{i+1} & ~~\textrm{if}~~ i = j \\
			t_jt_{j+1}s_{i-1} & ~~\textrm{if}~~ i = j+1 \\
			t_js_i & ~~\textrm{if}~~ i \ge j+2
				\end{array}\right.}
	\end{array}$$

\bigskip

In order to give an alternate presentation of the category $\fin^{op}$,
it is convenient to introduce the following operators first.

\pagebreak

\defn{
	The following maps in $\fin$ will be called the {\it standard cyclic
	permutations}. 
		$$z_i = z_i^{(n)}: [n] \lra [n] ~ \mfor n \ge 1 \mand 0 \le i \le n$$ 
		$$k \mapsto \left\{\begin{array} {cl}
				k+1 & \mfor 0 \le k \le i-1 \\
				0 & \mfor k = i \\
				k & \mfor k > i
			\end{array}\right.$$
	Note that $z_i$ is an $(i+1)$-cycle on the elements $0,1,\ldots, i$. In
	particular, $z_0$ is the identity. One may equivalently take the
	following formula in $\fin^{op}$ as a definition of the corresponding
	symmetric-simplicial operator $z_i$ for $n \ge 0 \mand 0 \le i \le n$.
		$$z_i = z_i^{(n)} := t_{i-1} \ldots t_1 t_0$$ 
}

\smallskip

The alternative presentation of $\fin^{op}$ given below keeps the
elementary face operators and transpositions as generators but
substitutes for the elementary degeneracies the following.

\defn{ \label{quasicodegens}
	The following maps in $\fin$ will be referred to as the {\it elementary
	quasi-codegeneracy} maps.
	$$u_i = u_i^{(n)}: [n+1] \lra [n] ~ \mfor n \ge 0 \mand 1 \le i \le n+1$$
	$$k \mapsto \left\{\begin{array} {cl}
				0 & \mfor k = 0 \textrm{~or~} i \\
				k & \mfor 1 \le k \le i-1 \\
				k-1 & \mfor k > i
			\end{array}\right.$$
	In particular, $u_1$ coincides with $s_0$. Note $u_0$ is not defined.
	One may equivalently define the {\it elementary quasi-degeneracy
	operators} $u_i$ $\fin^{op}$ in terms of the $s_i$ and $z_i$ by means of
	the following formula holding in $\fin^{op}$ for $i \ge 1$. 
		$$u_i := z_{i-1}^{-1}s_{i-1}z_{i-1}$$
}

\smallskip

The following theorem gives a presentation of $\fin^{op}$ in terms of
the generators $d_i, u_i, t_i$.

\pagebreak
\begin{thm}\label{QuasiIds}\label{altpresoffin}
	The generators $d_i$, $u_i$, and $t_i$ together with the 
	following relations constitute a presentation of $\fin^{op}$.
	$$\begin{array}{ll}
		{d_id_j = \left\{\begin{array}{cl}
				d_{j-1}d_{i} & ~~\mif~~ i < j \\
				d_jd_{i+1} & ~~\mif~~ i \ge j  
			\end{array}\right.} & \\[15pt]
		{d_iu_j = \left\{\begin{array}{cl}
				z_{j-1} & \phantom{0 \neq}~~~\mif~~ i=0 \\
				u_{j-1}d_{i} & ~~\mif~~          0 \neq  i < j \\
				\id & \phantom{0 \neq}~~~\mif~~  i = j \\
				u_jd_{i-1}   & \phantom{0 \neq}~~~\mif~~  i > j
			\end{array}\right.} &
		{u_id_j = \left\{\begin{array}{ll}
				d_{j+1}u_i & ~~\mif~~ i \le j \\
				d_ju_{i+1} & ~~\mif~~ i \ge j \neq 0 \\
				d_1u_{i+1}t_0 & ~~\mif~~ j = 0
			\end{array}\right.} \\[27pt]
		{u_iu_j = \left\{\begin{array}{cl}
				u_{j+1}u_{i} & ~~\mif~~ i \le j \\
				u_ju_{i-1} & ~~\mif~~ i > j
			\end{array}\right.} \\[13pt]
		{t_it_j = \left\{\begin{array}{cl}
			\id & ~~\mif~~ i = j \\
			t_jt_i & ~~\mif~~ \abs{i-j} \ge 2 \\
			(t_jt_i)^2 & ~~\mif~~ \abs{i-j} = 1
			\end{array}\right.} \\[20pt]
		{d_it_j = \left\{\begin{array}{cl}
			t_{j-1}d_{i} & ~~\mif~~ i < j \\
			d_{i+1} & ~~\mif~~ i = j \\
			d_{i-1} & ~~\mif~~ i = j+1 \\
			t_jd_i & ~~\mif~~ i \ge j+2
			\end{array}\right.}  &
		{t_id_j = \left\{\begin{array}{cl}
			d_jt_i & ~~\mif~~ i \le j-2 \\
			d_jt_{i+1}t_it_{i+1} & ~~\mif~~ i = j-1 \\
			d_jt_{i+1} & ~~\mif~~ i \ge j
			\end{array}\right.} \\[20pt]
		{t_iu_j = \left\{\begin{array}{cl}
				u_jt_i & ~~\mif~~ 0 \neq i \le j-2 \\
				u_{j-1} & ~~\mif~~ 0 \neq i = j-1 \\
				u_{j+1} & \phantom{0 \neq}~~~\mif~~ i = j \\
				u_jt_{i-1} & \phantom{0 \neq}~~~\mif~~ i > j \\
			\end{array}\right.} &
		{u_it_j = \left\{\begin{array}{cl}
				t_{j+1}u_i & ~~\mif~~ i \le j \\
				t_jt_{j+1}u_{i-1} & \begin{array}{r} ~\mif~~ i = j+1 \\ \mand j \neq 0\end{array} \\
				t_ju_i & \begin{array}{r} ~\mif~~ i \ge j+2 \\ \mand j \neq 0\end{array}
			\end{array}\right.} \\
		{t_0u_1 = u_1} \\
		{t_0u_it_0u_j = \left\{\begin{array}{rl}
			u_{j+1}t_0u_it_0  & \mif 2 \le i \le j \\
			u_jt_0u_{i-1}t_0 & \mif 2 \le j < i 
			\end{array}\right.} \\[7pt]
	\end{array}$$
\end{thm}

\bigskip

The next section is devoted to proving this theorem.

\smallskip

\rmk{
	For our purposes, this presentation has some advantages over that of
	Grandis. For instance, Corollary \ref{piugamma} is a consequence of the
	rule for $t_iu_j$ (contrast with the rule for $t_is_j$). The rule for
	$d_iu_j$ for $i > 0$ in particular is responsible for the vanishing of a
	great many brackets universally in symmetric hypercrossed complexes (this will
	be demonstrated in a forthcoming article).
\newline\indent
	The above presentation also has some notable disadvantages, particularly
	in the inability to move $t_0$ past any $u_i$ for $i \ge 2$, as well as
	in the identity $d_0u_i = z_{i-1}$, which makes the full definition
	of $d_0$ in symmetric hypercrossed complexes dependent on $t_0$.
}

\rmk{ 
	It is readily verified that all relations in the right column follow
	from the relations in the left column. All references to the statement of Theorem
	\ref{QuasiIds} will be understood as referring to relations of the
	left column only. 
}

\smallskip

Here are some other useful operators in $\fin^{op}$.

\defn{\label{replacementop}
	In the statement of Theorem \ref{QuasiIds}, note the overlapping conditions
	in the identities for $u_id_j$.  Indeed
	the equations
		$$u_id_i = d_{i+1}u_i = d_iu_{i+1} =: r_i$$
	hold for all $1 \le i \le n$.  We refer to the $r_i$ as {\it replacement operators}.
}

\begin{prp}
	For each $n \ge 1$, the replacement operators
		$$r_i: [n] \lra [n] \mfor 1\le i \le n$$
	constitute a family of mutually commuting idempotents in $\fin^{op}$.
		\begin{align*}
			r_i^2 &= r_i \\
			r_ir_j &= r_jr_i
		\end{align*}
\end{prp}
\prf{Proof.}{
	This is most easily verified using the following formula for $r_i$ as a
	function in $\fin$.
		$$r_i(k) = \left\{\begin{array}{cl}
					0 & \mif k = 0 \mathrm{~or~} i \\
					k & \mathrm{~otherwise~} 
				\end{array}\right.$$
	Alternatively, it is a fun exercise to prove the assertion using the 
	identities of Theorem \ref{QuasiIds} and Definition \ref{replacementop}.
}

\bigskip

%% file: SMCProof.tex

	\label{SMCProof}

In this section, we give the proof of Theorem \ref{QuasiIds}.

\bigskip

\prf{Proof of Theorem \ref{QuasiIds}.}{
	The framework for the proof is as follows.  Form the free category
	with objects $[n]$ for $n\ge 0$ and generators 
		$$\begin{array}{cll}
			d_i=d_i^{(n)}:&[n]\lra[n+1] 	&\mfor 0\le i\le n+1 \\
			t_i=t_i^{(n)}:&[n]\lra[n] 		&\mfor 0\le i\le n-1 \\
			s_i=s_i^{(n)}:&[n]\lra[n-1] 	&\mfor 0\le i\le n-1 \\
		\end{array}$$
	and let $\Q$ denote the quotient of this free category by those
	relations involving only the $d_i$ and $t_i$ (note these relations are
	common to both Grandis's presentation and the one proposed by 
	the theorem). According to \cite{Grandis}, the imposition on $\Q$ of
	the remaining relations of Grandis involving the $s_i$ produces the
	category $\fin^{op}$. Letting $u_i$ stand for
	$z_{i-1}^{-1}s_{i-1}z_{i-1}\in\Q$, we wish to show that the imposition
	of the relations in the statement of the theorem involving the $u_i$
	also produces $\fin^{op}$. For this it suffices to show that each
	relation involving the $u_i$ is a consequence of those involving the
	$s_i$ and those of $\Q$, and that each of Grandis's relations involving the $s_i$ is a
	consequence of those involving the $u_i$ and those of $\Q$. All the
	various statements constituting these assertions are proved in several
	propositions below, and the proof of the theorem is completed after
	that.
}

\begin{lem}\label{conjinsym}
	The following relations hold in the group $\sym[n]^{op}$.
	$$\begin{array}{ll}
		{z_i^{-1}z_j = \left\{\begin{array}{cl}
				z_jt_0z_{i+1}^{-1} & \mif i < j \\
				\id & \mif i=j \\
				z_{j+1}t_0z_{i}^{-1} & \mif i > j \\
			\end{array}\right.}
	\end{array} ~~~~
	z_{i}z_j = z_jz_{i+1}t_0 \mfor i < j$$
\end{lem}
\prf{Proof.}{
	Each of the above identities can be deduced from
		$$z_j^{-1}z_i z_j = z_{i+1}t_0 \mfor i < j$$
	so it suffices to prove the latter.  To see this, note that 
	since $i$ is less than $j$, the effect of the conjugation action of $z_j$ on 
		$$z_i = t_{i-1}\ldots t_0$$ 
	is to raise the subscript of each transposition by 1 (note we are working in
	$\sym[n]^{op}$ and not $\sym[n]$). The result is almost $z_{i+1}$, but
	$z_{i+1}$ has an extra $t_0$ at the end, so another $t_0$ is introduced
	to cancel it.
}

\begin{lem}\label{dzzds}
	The following identities hold in $\Q$ (hence also in $\fin^{op}$).
	$$\begin{array}{ll}
		{t_iz_{j} = \left\{\begin{array}{cl}
				z_{j}t_{i+1} & \mif i \le j-2 \\
				z_{j-1} & \mif i=j-1 \\
				z_{j+1} & \mif i = j \\
				z_{j}t_{i} & \mif i \ge j+1
			\end{array}\right.}
	\end{array} ~~~~
	\begin{array}{ll}
		{z_{i}^{-1}t_j = \left\{\begin{array}{cl}
				t_j z_i^{-1} & \mif i < j \\
				z_{i+1}^{-1} & \mif i = j \\
				z_{i-1}^{-1} & \mif i = j+1 \\
				t_{j+1}z_i^{-1} & \mif i \ge j+2 
			\end{array}\right.}
	\end{array}$$
		$$\begin{array}{ll}
			{d_i z_j = \left\{\begin{array}{cl}
				z_{j-1}d_{i+1} & \mif i < j \\
				d_0 & \mif i = j \\
				z_jd_i & \mif i > j
			\end{array}\right.} & ~~~~
			 {z_{i}^{-1}d_j  = \left\{\begin{array}{cl}
				d_j z_i^{-1} & \mif i < j \\
				d_{j+1}z_{i+1}^{-1} & \mif i \ge j
			\end{array}\right.}
		\end{array}$$
\end{lem}
\prf{Proof.}{
	For the upper batches of identities, use the relations of the Moore
	presentation in a straightforward manner. The proofs of the lower
	batches are similar to but easier than those of Lemma \ref{zuuzs} below.
	All are left to the reader.
}

\begin{prp}\label{Qfirst}
	In each of the following equivalences, imposing upon $\Q$ the relation
	on the left hand side produces the same result as imposing its
	correspondant on the right hand side.
	$$\begin{array}{lcl}
		d_is_j = s_{j-1}d_i \mfor i < j &\iff& d_iu_j = u_{j-1}d_i \mfor 0 \neq i \le j \\
		d_is_{i-1} = \id &\iff& d_iu_i = \id \\
		d_is_i = \id &\iff& d_0 u_{i+1} = z_i \\
		d_is_j = s_jd_{i-1} \mfor i \ge j+2 &\iff& d_iu_j = u_jd_{i-1} \mfor i > j \\
		\\
		t_is_j = s_jt_i \mfor i \le j-2 &\iff& t_i u_j = u_j t_i \mfor 0 \neq i \le j-2 \\
		t_{j-1}s_j = t_js_{j-1}t_{j-1} &\iff& t_iu_i = u_{i+1} \\
		t_{i+1}s_i = t_is_{i+1}t_i &\iff& t_iu_i = u_{i+1} \\
		t_is_j = s_jt_{i-1} \mfor i \ge j+2 &\iff& t_iu_j = u_jt_{i-1} \mfor i > j
%
	\end{array}	$$
\end{prp}
\prf{Proof.}{
	The verifications all follow the same pattern, so we prove the first one
	for illustration and leave the rest to the reader. For $i < j$, one has
	\begin{align*}
		d_i s_j &= d_i z_j u_{j+1} z_j^{-1} \tag{By Def. of $u_{j+1}$} \\
			&= z_{j-1} d_{i+1} u_{j+1} z_j^{-1} \tag{By Lemma \ref{dzzds}}
	\end{align*}
	and also 
	\begin{align*}
		s_{j-1}d_i 	&= z_{j-1} u_{j} z_{j-1}^{-1}d_i \tag{By Def. of $u_{j}$} \\
				&= z_{j-1} u_{j} d_{i+1} z_{j}^{-1} \tag{By Lemma \ref{dzzds}}
	\end{align*}
	and since the outermost terms of each of these results coincide and are
	invertible, one obtains (after a reparametrization) the following two-way
	implication as desired.
		$$d_is_j = s_{j-1}d_i \mfor i < j \iff d_iu_j = u_{j-1}d_i \mfor 0 \neq i \le j$$
}

\begin{lem}\label{zuuzs}
	Let $\Q^\p$ denote the category obtained by imposing on $\Q$ the 
	following relations from the statement of Theorem \ref{QuasiIds}.
		\begin{align}
			t_iu_j &= u_jt_i \mfor 0 \neq i \le j-2 \label{tuone}\\
			t_iu_{i+1} &= u_i \label{tutwo}\\
			t_iu_i &= u_{i+1} \mfor i\ge 1 \label{tuthr}\\
			t_iu_j &= u_jt_{i-1} \mfor i > j \label{tufou} \\
			t_0u_1 &= u_1 \label{tufiv}
		\end{align}
	Then the following identities hold in $\Q^\p$.
		$$\begin{array}{ll}
			{z_{j}^{-1} u_i = \left\{\begin{array}{cl}
				t_0 u_{i+1} t_0 z_{j-1}^{-1} & \mif i < j \\
				u_1 & \mif i = j \\
				t_0u_it_0z_{j}^{-1} & \mif i > j
			\end{array}\right.}& ~~~~
			{u_i z_j = \left\{\begin{array}{cl}
				z_{j+1} t_0 u_{i+1} t_0& \mif i \le j \\
				z_j t_0u_it_0 & \mif i > j
			\end{array}\right.}
		\end{array}$$
\end{lem}
\prf{Proof.}{
	The identities on the right can be directly deduced from those on the
	left, so we prove only the latter. The case $i<j$ is demonstrated as follows.
	\begin{align*}
		z_{j}^{-1} u_i  &= t_0 \ldots t_{j-1} u_i  \tag{Def. of $z_j$}\\
				&= t_0(t_1\ldots t_{i-1})t_i(t_{i+1}\ldots t_{j-1}) u_i \\
				&= t_0(t_1\ldots t_{i-1})t_i u_i (t_{i}\ldots t_{j-2}) \tag{By \eqref{tufou}} \\
				&= t_0(t_1\ldots t_{i-1}) u_{i+1} (t_{i}\ldots t_{j-2}) \tag{By \eqref{tuthr}} \\
				&= t_0 u_{i+1} (t_1\ldots t_{i-1})(t_{i}\ldots t_{j-2}) \tag{By \eqref{tuone}} \\
				&= t_0 u_{i+1} t_0t_0(t_1\ldots t_{i-1})(t_{i}\ldots t_{j-2}) \tag{$t_0^2 = \id$} \\
				&= t_0 u_{i+1} t_0z_{j-1}^{-1} \tag{Def. of $z_{j-1}$}
	\end{align*}
	The proof of the case $i > j$ is similar but easier. The case $i=j$
	follows from \eqref{tufiv} and repeated application of \eqref{tutwo}.
}

\begin{lem}\label{uus}
	Let $\Q^\p$ be as in Lemma \ref{zuuzs}.
	Then the following holds in $\Q^\p$.
		$${u_iu_j = \left\{\begin{array}{cl}
				u_{j+1}u_{i} & ~~\mif~~ i \le j \\
				u_ju_{i-1} & ~~\mif~~ i > j
			\end{array}\right.}$$
\end{lem}
\prf{Proof.}{
	Calculate as follows for $i\le j$.
	\begin{align*}
		u_iu_j  &= z_i u_1 z_j u_1 \tag{By Lemma \ref{zuuzs}}\\
			&= z_i z_{j+1} t_0 u_1 t_0 u_1 \tag{By Lemma \ref{zuuzs}}\\
			&= z_{j+1} z_{i+1}  t_0 u_1 t_0 u_1 \tag{By Lemma \ref{conjinsym}}\\
			&= z_{j+1} u_1 z_i u_1 \tag{By Lemma \ref{zuuzs}}\\
			&= u_{j+1} u_{i} \tag{By Lemma \ref{zuuzs}}
	\end{align*}
	The case $i>j$ follows immediately from the case $i\le j$.
}

\begin{prp}\label{Qsecond}
	Let $\Q^\p$ be as in Lemma \ref{zuuzs}.
	Then the following holds in $\Q^\p$.
		$$t_is_i = s_i \mforall i$$
\end{prp}
\prf{Proof.}{
	One computes as follows.
	\begin{align*}
		t_is_i  &= t_i z_iu_{i+1}z_i^{-1} \tag{Def. of $u_{i+1}$}\\
			&= z_{i+1}u_{i+1}z_i^{-1} \tag{By Lemma \ref{conjinsym}} \\
			&= z_{i+1}z_{i+1}u_{1}z_i^{-1} \tag{By Lemma \ref{zuuzs}} \\
			&= z_{i}z_{i+1}t_0u_{1}z_i^{-1} \tag{By Lemma \ref{conjinsym} with $j=i+1$} \\
			&= z_{i}z_{i+1}u_{1}z_i^{-1} \tag{$t_0u_1 = u_1$} \\
			&= z_{i}u_{i+1}z_i^{-1} \tag{By Lemma \ref{zuuzs}} \\
			&= s_i \tag{Def. of $u_{i+1}$}
	\end{align*}
}

\begin{prp}\label{sstutu}
	Let $\Q^\p$ be as in Lemma \ref{zuuzs}.  Then the following equivalences
	of algebraic relations hold (in the same sense as in Lemma \ref{Qfirst}).
	$$\begin{array}{lcl}
		s_is_j = s_{j+1}s_i \mfor i < j & \iff & t_0u_it_0u_j = u_{j+1}t_0u_it_0 \mif 2 \le i \le j\\
		s_is_i = s_{i+1}s_i \mforall i \ge 0 & \iff & u_1u_{j} = u_{j+1}u_1 \mforall j\ge1
	\end{array}$$
\end{prp}
\prf{Proof.}{
	Calculate as follows for $i<j$.
	\begin{align*}
		s_is_j  &= z_i u_{i+1} z_i^{-1} z_j u_{j+1} z_j^{-1} \tag{By Def. \ref{quasicodegens}} \\
			&= z_i u_{i+1} z_j t_0 z_{i+1}^{-1} u_{j+1} z_j^{-1} \tag{By Lemma \ref{conjinsym}} \\
			&= z_i z_{j+1} t_0 u_{i+2} t_0t_0t_0 u_{j+1} t_0 z_{i+1}^{-1} z_j^{-1} \tag{By Lemma \ref{zuuzs} twice} \\
			&= (z_i z_{j+1} t_0) (u_{i+2} t_0 u_{j+1} t_0) (z_{i+1}^{-1} z_j^{-1})
	\end{align*}
	\begin{align*}
		s_{j+1}s_i &= z_{j+1}u_{j+2}z_{j+1}^{-1}z_i u_{i+1} z_i^{-1}\tag{By Def. \ref{quasicodegens}} \\
				&= z_{j+1}u_{j+2}z_{i+1}t_0z_{j+1}^{-1} u_{i+1} z_i^{-1}\tag{By Lemma \ref{conjinsym}} \\
				&= z_{j+1}z_{i+1}t_0u_{j+2}t_0t_0t_0 u_{i+2} t_0z_{j}^{-1}z_i^{-1}  \tag{By Lemma \ref{zuuzs} twice} \\
				&= (z_{j+1}z_{i+1})(t_0u_{j+2}t_0 u_{i+2}) (t_0z_{j}^{-1}z_i^{-1})
	\end{align*}
	Now by Lemma \ref{conjinsym}, the respective outer terms of the two
	expressions coincide, and since these terms are invertible, one obtains
	the equality of the inner terms, that is
		$$s_is_j = s_{j+1}s_i \mfor i < j \iff t_0u_it_0u_j = u_{j+1}t_0u_it_0 \mfor 2 \le i \le j$$
	Similarly, for the case $i=j$ one computes
	\begin{align*}
		s_is_i  &= z_i u_{i+1} z_i^{-1} z_i u_{i+1} z_i^{-1} \tag{By Def. \ref{quasicodegens}} \\
			&= z_i z_{i+1} u_1 u_{i+1} z_i^{-1} \tag{By Lemma \ref{zuuzs}}
	\end{align*}
	\begin{align*}
		s_{i+1}s_i 	&= z_{i+1}u_{i+2}z_{i+1}^{-1}z_iu_{i+1}z_i^{-1} \tag{By Def. \ref{quasicodegens}}\\
				&= z_{i+1}u_{i+2}z_{i+1}t_0z_{i+1}^{-1}u_{i+1}z_i^{-1} \tag{By Lemma \ref{conjinsym}}\\
				&= z_{i+1}z_{i+1}t_0u_{i+2}t_0t_0z_{i+1}^{-1}u_{i+1}z_i^{-1} \tag{By Lemma \ref{zuuzs}}\\
				&= z_{i+1}z_{i+1}t_0u_{i+2}t_0t_0u_{1}z_i^{-1} \tag{By Lemma \ref{zuuzs}}\\
				&= z_{i+1}z_{i+1}t_0u_{i+2}u_{1}z_i^{-1} \tag{By Lemma \ref{zuuzs}}
	\end{align*}
	and similarly as before one obtains the following.
		$$s_is_i = s_{i+1}s_i \mforall i \ge 0 \iff u_1u_{j} = u_{j+1}u_1 \mforall j\ge1$$
}

In order to finish the
proof of Theorem \ref{QuasiIds}, it is convenient to introduce the following notation.

\defn{ \label{Rdu}
	Let the symbol ${\bf R}[d,u]$ denote the set of
	relations stated in Theorem \ref{QuasiIds} involving only the operators $d_i$ and $u_j$.
	Similarly use the symbols ${\bf R}[u,u]$ and ${\bf R}[t,u]$. It is
	also convenient to write ${\bf R}[t_+,u]$ for those relations of ${\bf
	R}[t,u]$ involving only $t_i$ with $i>0$ and ${\bf R}[t_0,u]$ for those
	relations of ${\bf R}[t,u]$ involving $t_0$ but not $t_i$ with $i>0$.
\newline\indent
	Additionally ${\bf R}[d,s]$, ${\bf R}[s,s]$ and ${\bf R}[t,s]$ are used to
	refer to the analogous sets of relations from Grandis's presentation.
}

\prf{Proof of Theorem \ref{QuasiIds} (continued).}{
	For one direction, assume that all of Grandis's relations hold, so that
	we are working in the category $\fin^{op}$. 

	From Proposition \ref{Qfirst} all identities ${\bf R}[d,u]$ and ${\bf
	R}[t_+,u]$ are obtained. These identities fulfill the hypotheses of
	Lemmas \ref{zuuzs} and \ref{uus}, and thus ${\bf R}[u,u]$ is obtained.
	Finally, these same lemmas enable us to apply Proposition \ref{sstutu}
	and so ${\bf R}[t_0,u]$ is also obtained (with the exception of 
	$t_0u_1 = u_1$, but this is just the same as $t_0s_0 = s_0$).

	For the reverse direction, assume all identities ${\bf R}[d,u]$, ${\bf
	R}[t,u]$ and ${\bf R}[u,u]$ are imposed on $\Q$. In the resulting
	category, all identities ${\bf R}[d,s]$ and ${\bf R}[t,s]$ obtain by
	Propositions \ref{Qfirst} and \ref{Qsecond}. By Proposition
	\ref{sstutu}, the identities ${\bf R}[s,s]$ also obtain.
}

\rmk{
	It is a corollary of Propositions \ref{Qfirst} and \ref{Qsecond} that,
	in Grandis's presentation \cite{Grandis}, the relations $t_is_i = s_i$
	for $i > 0$ are redundant. Tracing this, one finds that they are a
	consequence of the relation $t_0s_0 = s_0$ as well as the other
	relations for exchanging $t_i$ and $s_j$.
\newline\indent
	Similarly and perhaps surprisingly, 
	Lemma \ref{uus} says that the relations for $u_iu_j$ are
	redundant in the alternate presentation of Theorem \ref{QuasiIds}. They
	are a consequence of the relations for exchanging $t_i$ and $u_j$ as
	well as of the Moore relations for the $t_i$.
}

%% file: SymAlgMultIndex.tex
	\label{SymAlgMultIndex}

In the subsections of this final section we collect together a number of facts about the
algebraic structure of the category $\fin$, viewed from the point of
view of the alternative presentation given in Theorem
\ref{altpresoffin}.  For this purpose we introduce the following notational device.

\defn{
A {\it multi-index $\alpha$ of length $k$ and dimension $\le n$} is a strictly increasing sequence
of indices
	$$\alpha = \big\{i_1 < \ldots < i_k \big\}$$
satisfying $1 \leq i_p \leq n$ for all $p$. The length
$\vert\alpha\vert$ of $\alpha$ is the number $k$ of indices in $\alpha$.
}
\defn{
The {\it quasi-codegeneracy} $u_\alpha\in\fin$
corresponding to $\alpha$ is the composition of elementary
quasi-codegeneracies
	$$u_\alpha : [n] \longrightarrow [n-\vert\alpha\vert]$$
	$$u_\alpha := u_{i_1}u_{i_2}\ldots u_{i_k}$$
(note the indices increase from left to right) 
and the {\it coface} $d_\alpha\in\fin$ corresponding to $\alpha$ is
	$$d_\alpha : [n-\vert\alpha\vert] \longrightarrow [n]$$
	$$d_\alpha := d_{i_k}d_{i_{k-1}}\ldots d_{i_1}$$
(note the indices decrease from left to right).
}

\subsection{The Quasi-Monotonic Functions}

The goal of this subsection is to characterize the functions in $\fin$
obtained as compositions of the form $d_\alpha u_\beta$ as the {\it
quasi-monotonic functions} (defined below), to show that they constitute a
subcategory of $\fin$ and finally to prove that the expressions $d_\alpha
u_\beta$ themselves constitute a family of unique factorizations for
that subcategory. This amounts to an analog of the usual unique
factorization theorem for $\ord$ with the $u_i$ taking the place of the
$s_i$, we which recall here (see \cite{May}, \cite{MacLane},
\cite{Lamotke} for a proof). 

\begin{prp} \label{standuniqfact}
	Each simplicial operator $f:[n]\ra[m]\in \ord$ has a unique factorization
		$$f = d_{i_l}\ldots d_{i_1} s_{j_1} \ldots s_{j_k}$$
	with strictly increasing sequences of indices
		$$0 \le i_1 < \ldots < i_k \le m$$
		$$0 \le j_1 < \ldots < j_l \le n,$$
	that is, reading from left to right, the indices of the degeneracies decrease and
	the indices of the faces increase.  
\end{prp}

\bigskip

The following property will ultimately be shown to characterize
functions of the form $d_\alpha u_\beta$. 

\defn{
	Let a function $f\in\fin$ be called {\it quasi-monotonic} if it 
	satisfies the following two conditions.
	\begin{enumerate}
		\item[{QM1.}]{$f(0) = 0$}
		\item[{QM2.}]{$f(p) \neq 0 \neq f(q) \text{~and~} p < q \implies f(p) < f(q)$}
	\end{enumerate}
	The second condition says that
	$f$ is strictly increasing outside of $f^{-1}(0)$.
}

\smallskip

\begin{prp}\label{quasimonoclosedcomp}
	The quasi-monotonic functions are closed under composition.
\end{prp}
\prf{Proof.}{
	Let $f$ and $g$ be quasi-monotonic.  It suffices to check for $p<q$ that 
		$$(g\circ f)(p) \neq 0 \neq (g\circ f)(q) \implies (g\circ f)(p) < (g\circ f)(q).$$
	By QM1 for $g$ it must be the case that $f(p) \neq 0 \neq f(q)$,
	and it follows that $f(p) < f(q)$ since $f$ is quasi-monotonic.  Then condition 
	QM2 for $g$ applies to give exactly $g(f(p)) < g(f(q))$.
}

\begin{lem}\label{ualphavalues}
	Let $\alpha$ be a multi-index and let 
	the nonnegative integer $p$ belong to the domain of $u_\alpha$.
	Then $u_\alpha(p) = 0$ if and only if $p = 0$ or $p$ belongs to $\alpha$.
	Moreover, if $p$ is neither $0$ nor in $\alpha$ then 
		$$u_\alpha(p) = p - \#\big\{~ i\in\alpha ~\big\vert~ i<p ~\big\}.$$
\end{lem}
\prf{Proof.}{
	We start with the ``if'' direction.  First, if $p$ is $0$, then $u_\alpha(p) = 0$ since by
	Definition \ref{quasicodegens}, all quasi-codegeneracies send 0 to 0. Now assume $p$ belongs to $\alpha$.
	Then one may factor $u_\alpha$ as $u_{\alpha_{<p}}u_p u_{\alpha_{>p}}$ where $\alpha_{<p}$ consists of those
	indices in $\alpha$ that are less than $p$ and $\alpha_{>p}$ consists of those
	indices in $\alpha$ that are greater than $p$.  By Definition \ref{quasicodegens}, all $u_i$ with $i>p$ send
	$p$ to $p$, so that $u_{\alpha_{>p}}(p) = p$.  Then we have 
		\begin{align*}
			u_\alpha(p) &= u_{\alpha_{<p}}u_p u_{\alpha_{>p}}(p) \\
			&= u_{\alpha_{<p}}u_p(p) \\
			&= u_{\alpha_{<p}}(0) \\
			&= 0 
		\end{align*}
	as claimed.

	We turn to prove the ``only if'' direction, so we assume $p$ is neither $0$ nor in $\alpha$.
	Then one may factor $u_\alpha$ as $u_{\alpha_{<p}}u_{\alpha_{>p}}$ with $\alpha_{<p}$ and $\alpha_{>p}$ 
	having the same meanings as above.  As before $u_{\alpha_{>p}}(p) = p$, 
	so that $u_\alpha(p) = u_{\alpha_{<p}}(p)$.  If $\alpha_{<p}$ is empty then $u_{\alpha_{<p}}$
	is the identity and we conclude $u_\alpha(p) = p \neq 0$ as required.
	If $\alpha_{<p}$ is not empty, say $\alpha_{<p} = \{i_1,i_2,\ldots,i_j\}$.  Then by 
	Definition \ref{quasicodegens} we can evaluate 
		\begin{align*}
			u_{\alpha_{<p}}(p) &= u_{i_1}u_{i_2}\ldots u_{i_j}(p) \\
			&= u_{i_1}u_{i_2}\ldots u_{i_{j-1}}(p-1) \\
			&\vdots \\
			&= u_{i_1}u_{i_2}\ldots u_{i_{j-l}}(p-l) \\
			&\vdots \\
			&= p-j.
		\end{align*}
	where one notes that at each step, the argument $p-l$ decreases by 1 while the rightmost index $i_{j-l}$
	decreases by at least 1, so that $i_{j-l} < p-l$ holds for all $l$ and therefore the calculation
	may always proceed to the next step.
	Since $j$ is the number of indices in $\alpha_{<p}$ that are less than $p$, and since all
	indices in $\alpha_{<p}$ are between 1 and $p-1$ inclusive, we deduce that $j$ is at most $p-1$.
	Then $u_\alpha(p) = p-j \ge 1$ and we conclude that $u_\alpha(p)$ is not $0$ as claimed.

	The final assertion of the Lemma can be read out of the proof of the ``only if'' direction just given.
}

\begin{lem}\label{qmsurjunique}
	Quasi-monotonic surjections are uniquely determined by their zeros.
\end{lem}
\prf{Proof.}{
	Let $h:[n]\lra[m]$ be a quasi-monotonic surjection.  Note $h$ restricts
	to a surjection $[n]\setminus h^{-1}(0) \lra [m] \setminus \{0\}$ and by QM2 
	the restriction must be strictly monotonic, therefore also a bijection.  
	Hence $h$ is uniquely determined on $[n]\setminus h^{-1}(0)$
	and therefore on all of $[n]$.  
}

\begin{prp}\label{charquassurjs}
	The quasi-monotonic injections are precisely the functions $d_\alpha$.  The quasi-monotonic
	surjections are precisely the functions $u_\alpha$.
\end{prp}
\prf{Proof.}{
	That all functions of the form $d_\alpha$ or $u_\beta$ are quasi-monotonic
	follows from Lemma \ref{quasimonoclosedcomp} and the fact that, with the
	exception of $d_0$, all $d_i$ and $u_i$ are quasi-monotonic. Moreover $d_0$ is
	not a factor of $d_\alpha$ since multi-indices $\alpha$ by definition do
	not contain 0.

	Let $g$ be a quasi-monotonic injection.  Then $g$ must be monotonic because 
	by QM2, $g$ is increasing off of $g^{-1}(0)$ and by QM1 $g^{-1}(0)$ is $\{0\}$, i.e., $g$ is increasing 
	on the rest of the domain of $g$.  Then $g$ has a unique factorization 
	$g = d_{i_k} \ldots d_{i_1}$ in $\ord$ with $i_1 < \ldots < i_k$, and moreover $i_1$
	cannot be 0 because then $0$ would not be in the image of $g$.  So $\alpha = \{ i_1 < \ldots < i_k \}$
	satisfies the definition of multi-index and $g$ is equal to $d_\alpha$.

	Finally let $h$ be a quasi-monotonic surjection.  
	Let $\alpha$ be the multi-index consisting of the zeros of $h$, excluding 0 itself.
	By Lemma \ref{ualphavalues}, the function $u_\alpha$ has precisely the same zeros as $h$.
	Then by Lemma \ref{qmsurjunique}, $h$ must coincide with $u_\alpha$.
}

\begin{prp}\label{charquasfuncs}
	The quasi-monotonic functions are precisely the functions of the form $d_\alpha u_\beta$.  
\end{prp}
\prf{Proof.}{
	That all functions of the form $d_\alpha u_\beta$ are quasi-monotonic
	follows from Lemma \ref{quasimonoclosedcomp} just as in the preceding proposition.

	To see that every quasi-monotonic $f$ has the form $d_\alpha u_\beta$,
	first factor $f$ as $f = g\circ h$ where $g$ is a monotonic injection
	and $h$ is a surjection (this is possible for any $f$). Then $g$ has a
	factorization $d_{i_k} \ldots d_{i_1}$ by Proposition
	\ref{standuniqfact}. Since 0 is in the image of $f$, $d_0$ does not
	occur, that is, none of the $i_j$ is 0, and hence we may write down the quasi-monotonic
	surjection $u_{i_1}\ldots u_{i_k}$. It is a left-inverse for $g$,
	so by composing with it one discovers that 
		$$h = (u_{i_1}\ldots u_{i_k})\circ g \circ h = (u_{i_1}\ldots u_{i_k})\circ f$$ 
	is quasi-monotonic by Lemma \ref{quasimonoclosedcomp}. Then $h$ is a quasi-monotonic surjection,
	so by Proposition \ref{charquassurjs}, $h$ has the form
	$u_{j_1}\ldots u_{j_l}$.
}

\begin{prp}\label{charquasfuncs}
	Quasi-monotonic functions $f\in \fin$ have unique factorizations
	of the form $f = d_\alpha u_\beta$.
\end{prp}
\prf{Proof.}{
	Existence of factorizations was demonstrated in the previous proposition.
	To prove uniqueness, let a quasi-monotonic function $f$ have two factorizations
		$$f = d_\alpha u_\beta  = d_{\alpha^\p} u_{\beta^\p}$$
	where $\alpha,\alpha^\p,\beta,\beta^\p$ are multi-indices.
	It follows that
		$$\mathrm{Im}(d_\alpha) = \mathrm{Im}(d_\alpha u_\beta) = 
			\mathrm{Im}(d_{\alpha^\p} u_{\beta^\p}) = \mathrm{Im}(d_{\alpha^\p})$$
	and since monotonic injections are uniquely determined by their images, one concludes
	$d_\alpha = d_{\alpha^\p}$ and then $\alpha = \alpha^\p$ by Proposition \ref{standuniqfact}.  
	Then cancelling these by composing both sides with a left-inverse,
	one obtains $u_\beta = u_{\beta^\p}$ and finally $\beta = \beta^\p$ by Lemmas \ref{ualphavalues}
	and \ref{qmsurjunique}.
}

\subsection{Exchanging Transpositions and Quasidegeneracies}

In this subsection we consider interchanging permutations with a
quasidegeneracy $u_\gamma$. 
The reader is warned that we state results for $\fin^{op}$ instead of for $\fin$.
\smallskip

The following definition is motivated by the
effect of repeatedly using the identity for $t_iu_j$ from Theorem
\ref{QuasiIds} to push $t_i$ to the right across the factors of
$u_\gamma$ one at a time.

\defn{\label{multitransop}
	For any multi-index $\gamma$ and index $i\neq 0$, let $t_{i\sim\gamma}$
	stand for the following permutation in $\fin^{op}$.
		$$t_{i\sim\gamma} := \left\{\begin{array}{cl}
				\id & \mif i \text{~or~} i+1 \in \gamma \\
				t_{i^\p} \mfor i^\p := i - \#\setof{j\in\gamma}{j<i} & 
						\mif i\notin\gamma \mand i+1\notin\gamma
			\end{array}\right.$$
}


\bigskip

In the next lemma, $t_i(\gamma)$ denotes the multi-index obtained by applying $t_i$
as a function to the elements of $\gamma$.  Note that the indices in $t_i(\gamma)$ are
to be rearranged in sequential order even if $t_i$ puts the indices of $\gamma$ out of order.

\begin{lem}
	The following identity holds in $\fin^{op}$ for any  multi-index $\gamma$ and
	index $i \neq 0$.
		$$t_iu_\gamma = u_{t_i(\gamma)} t_{i\sim\gamma}$$
\end{lem}
\prf{Proof.}{
	Write 
		$$u_\gamma = u_\zeta u_\theta u_\alpha$$
	where $\zeta$ consists of the indices of $\gamma$ greater than or equal to $i+2$, $\theta$ consists
	of those of $i$ and $i+1$ belonging to $\gamma$, and $\alpha$ consists of those
	indices of $\gamma$ less than $i$.  

	From the following identities from Theorem \ref{QuasiIds}
		\begin{align*}
			t_i u_i &= u_{i+1} \\
			t_i u_{i+1} &= u_i \\
			t_i u_{i+1} u_i &= u_i u_i = u_{i+1} u_i
		\end{align*}
	one concludes
		$$t_i u_\theta = u_{t_i(\theta)}$$
	whenever $\theta$ is not empty.  In this case, one calculates
		\begin{align*}
			t_i u_\gamma &= t_i u_\zeta u_\theta u_\alpha \\
			&= u_\zeta t_i u_\theta u_\alpha \tag{By Theorem \ref{QuasiIds}} \\
			&= u_\zeta u_{t_i(\theta)} u_\alpha \\
			&= u_{t_i(\gamma)}
		\end{align*}
	as desired.

	If $\theta$ is empty, that is, neither $i$ nor $i+1$ belongs to $\gamma$, 
	then the calculation becomes
		\begin{align*}
			t_i u_\gamma &= t_i u_\zeta u_\alpha \\
			&= u_\zeta t_i u_\alpha \tag{By Theorem \ref{QuasiIds}} \\
			&= u_\zeta u_\alpha t_{i-\abs{\alpha}} \tag{By Theorem \ref{QuasiIds}}\\
			&= u_\gamma t_{i-\abs{\alpha}} \\
			&= u_{t_i(\gamma)} t_{i\sim\gamma}
		\end{align*}
	also as desired.
}

\bigskip

For the final corollary below, we let $\sym_n^\p$ denote the group of permutations of $[n]$
leaving $0\in[n]$ fixed.

\begin{cor} \label{piugamma}
	For any permutation $\pi\in(\sym_n^\p)^{op}$ and multi-index $\gamma$ 
	of length $k$, there exists a permutation $\pi^\p \in (\sym_{n-k}^\p)^{op}$ such that
		$$\pi u_\gamma = u_{\pi^{-1}(\gamma)} \pi^\p$$
\end{cor}
\prf{Proof.}{
	Factor $\pi$ as a product of operators $t_i$ and use the previous lemma
	to push each one past $u_\gamma$. In this process, all
	transpositions in the factorization of $\pi$ pile up in order in the
	subscript of $u_\gamma$. Since $\pi$ belongs to $\fin^{op}$, its
	factorization occurs in the order opposite to that of its factorization in the permutation
	group $\sym_n^\p$. Applying the transpositions as functions in this
	reversed order to $\gamma$ is therefore the same as applying the permutation
	$\pi^{-1}$ to $\gamma$.	
}

%% file: Preprint_4_QuasiMonoticPresFin.bbl
\providecommand{\bysame}{\leavevmode\hbox to3em{\hrulefill}\thinspace}
\providecommand{\MR}{\relax\ifhmode\unskip\space\fi MR }
\providecommand{\MRhref}[2]{%
  \href{http://www.ams.org/mathscinet-getitem?mr=#1}{#2}
}
\providecommand{\href}[2]{#2}
\begin{thebibliography}{DHK85}

\bibitem[Ant10]{preprint2}
Eric~R. Antokoletz, \emph{Nonabelian dold-kan decompositions for simplicial and
  symmetric-simplicial groups}, available online at http://www.arxiv.org
  (2010).

\bibitem[CC91]{CarrascoCegarra}
P.~Carrasco and A.~M. Cegarra, \emph{Group-theoretic algebraic models for
  homotopy types}, J. Pure Appl. Algebra \textbf{75} (1991), no.~3, 195--235.
  \MR{MR1137837 (93b:55026)}

\bibitem[DHK85]{DwyerHopkinsKan}
W.G. Dwyer, M.~J. Hopkins, and D.~M. Kan, \emph{The homotopy theory of cylic
  sets}, Trans. Amer. Math. Soc. \textbf{291} (1985), no.~1, 281--289.
  \MR{MR770723 (86m:55014)}

\bibitem[GJ99]{GoerssJardine}
Paul~G. Goerss and John~F. Jardine, \emph{Simplicial homotopy theory}, Progress
  in Mathematics, vol. 174, Birkh\"auser Verlag, Basel, 1999. \MR{MR1711612
  (2001d:55012)}

\bibitem[Gra88]{Grandis2}
Marco Grandis, \emph{Simplicial toposes and combinatorial homotopy}, Dip. Mat.
  Univ. Genova, Preprint 400 (1999)., 1988.

\bibitem[Gra01a]{Grandis}
\bysame, \emph{Finite sets and symmetric simplicial sets}, Theory Appl. Categ.
  \textbf{8} (2001), 244--252 (electronic). \MR{MR1825431 (2002c:18010)}

\bibitem[Gra01b]{Grandis3}
\bysame, \emph{Higher fundamental functors for simplicial sets}, Cahiers
  Topologie G\'eom. Diff\'erentielle Cat\'eg. \textbf{42} (2001), no.~2,
  101--136. \MR{MR1839359 (2002f:18026)}

\bibitem[Gra02]{Grandis1}
\bysame, \emph{An intrinsic homotopy theory for simplicial complexes, with
  applications to image analysis}, Appl. Categ. Structures \textbf{10} (2002),
  no.~2, 99--155. \MR{MR1891107 (2003a:18014)}

\bibitem[Gra03]{GrandisNew}
\bysame, \emph{Higher fundamental groupoids for spaces}, Topology Appl.
  \textbf{129} (2003), no.~3, 281--299. \MR{MR1962985 (2004c:55036)}

\bibitem[Gro83]{Groth}
Alexander Grothendieck, \emph{Pursuing stacks}, 1983.

\bibitem[Lam68]{Lamotke}
Klaus Lamotke, \emph{Semisimpliziale algebraische {T}opologie}, Die Grundlehren
  der mathematischen Wissenschaften, Band 147, Springer-Verlag, Berlin, 1968.
  \MR{MR0245005 (39 \#6318)}

\bibitem[Law88]{Lawvere}
F.W. Lawvere, \emph{Toposes generated by codiscrete objects, in combinatorial
  topology and functional analysis}, Notes for Colloquium lectures given at
  North Ryde, NSW, Aus (1988), at Madison, Wisconsin, USA (1989), and at
  Seminario Matematico e Fisico, Milano, Italy (1992)., 1988.

\bibitem[May67]{May}
J.~Peter May, \emph{Simplicial objects in algebraic topology}, Van Nostrand
  Mathematical Studies, No. 11, D. Van Nostrand Co., Inc., Princeton,
  N.J.-Toronto, Ont.-London, 1967. \MR{MR0222892 (36 \#5942)}

\bibitem[ML70]{MacLane}
Saunders Mac~Lane, \emph{The {M}ilgram bar construction as a tensor product of
  functors}, The Steenrod Algebra and its Applications (Proc. Conf. to
  Celebrate N. E. Steenrod's Sixtieth Birthday, Battelle Memorial Inst.,
  Columbus, Ohio,1970), Lecture Notes in Mathematics, Vol. 168, Springer,
  Berlin, 1970, pp.~135--152. \MR{MR0273618 (42 \#8495)}

\end{thebibliography}
